\documentclass[12pt]{article} %arxiv

\usepackage{amsfonts}%arxiv

\newcommand{\C}{\mathbb{C}}

\newcommand{\Q}{\mathbb{Q}}

\newcommand{\Z}{\mathbb{Z}}

\newtheorem{thm}{Theorem}

\newcommand{\dd}{\;\mathrm{d}}%differential
\newcommand{\e}{\mathrm{e}}%Exponential

\newcommand{\qed}{$\Box$}

\newcommand{\ii}{\mathrm{i}} %Imaginary unit
\begin{document}

\title{On a conjecture by Boyd} %arxiv

\author{Matilde N. Lal\'in\footnote{Department of Mathematical and Statistical Sciences, University of Alberta, 632 Central Academic Building,Edmonton, AB T6G 2G1, Canada
\texttt{mlalin@math.ualberta.ca}}} %arxiv

%\markboth{Matilde N. Lal\'in} %journal
%{On a conjecture by Boyd} %journal

%\catchline{}{}{}{}{}%journal
%
%%%%%%%%%%%%%%%%%%%%%%%%%%%%%%%%%%%%%%%%%%%%%%%%%%%%%%%%%%%%%%%%%%%%

%\title{ON A CONJECTURE BY BOYD} %journal

%\author{MATILDE N. LAL\'IN} %journal

%\address{ Department of Mathematical and Statistical Sciences, University of Alberta, 632 Central Academic Building,Edmonton, AB T6G 2G1, Canada\\
%\email{mlalin@math.ualberta.ca}} %journal

\maketitle

%\begin{history}%journal
%\received{(Day Month Year)}
%\accepted{(Day Month Year)}
%\comby{xxx}
%\end{history}

\begin{abstract} The aim of this note is to prove the Mahler measure identity $m(x+x^{-1}+y+y^{-1}+5) = 6 m(x+x^{-1}+y+y^{-1}+1)$ which was conjectured by Boyd. The proof is achieved by proving relationships between regulators of both curves.
\end{abstract}

%\keywords{ Mahler measure, elliptic curves, elliptic dilogarithm, regulator}%journal

\noindent{\it keywords}: Mahler measure, elliptic curves, elliptic dilogarithm, regulator %arxiv
\noindent{\it 2000 Mathematics Subject Classification}: 2000: 11R09, 19F27 %arxiv
%\ccode{Mathematics Subject Classification 2000: 11R09, 19F27}%journal

\section{Introduction}

Boyd \cite{Bo} studied the Mahler measure of families of polynomials. In particular, he considered the two-variable family
\[P_k(x,y) = x +\frac{1}{x} + y +\frac{1}{y} + k.\]
The zeros of $P_k(x,y)$ correspond, generically to a curve of genus 1. Let $E_k$ denote the elliptic curve corresponding to the algebraic closure of $P_k(x,y)=0$.

Recall that the (logarithmic) Mahler measure of a non-zero Laurent polynomial, $P(x_1, \dots, x_n)$, with complex coefficients is defined as
\[m(P)=\int_0^1 \dots \int_0^1 \log \left|P\left(\e^{2\pi\ii t_1}, \dots, \e^{2\pi\ii t_n}\right)\right| \dd t_1 \dots t_n.\]

Let us denote $m(k) :=m(P_k)$. Boyd computed $m(k)$ for $k$ a positive integer less than or equal to 100 (it is easy to see that the Mahler measure does not depend on the sign of $k$ for this family). He found that 
\begin{equation} \label{eq:Boyd}
m(k) \stackrel{?}{=} r_kL'(E_k,0),
\end{equation}
where $r_k$ is a rational number and the question mark stands for an equality that has only been stablished numerically (typically to at least 50 decimal places). 

The case with $k=1$ (resulting in $r_k=1$) was considered in detail by Deninger \cite{D}, who found an explanation for such a formula by relating it to evaluations of regulators in the context of the Bloch--Beilinson conjectures. Rodriguez-Villegas \cite{RV} also considered this family in the context of the Bloch-Beilinson conjectures, including more general cases where $k^2 \in \Q$. He was able to prove identities for the cases where the Bloch--Beilinson conjectures are known to be true, such as when $E_k$ has complex multiplication.

When the curves $E_{k_1}$ and $E_{k_2}$ are isogenous, their $L$-functions coincide. One can then compare the values in equation $(\ref{eq:Boyd})$ and conjecture identities of the form $r_{k_2} m(k_1) = r_{k_1} m(k_2)$. For example, 
\begin{thm}
\begin{equation} \label{eqA}
m(8) = 4 m(2),
\end{equation}
\begin{equation} \label{eqB}
m(5)=6m(1).
\end{equation}
\end{thm}

The first identity was proved in \cite{LR}. In this note, we prove the second one. 

\section{Functional Identities}

Functional identities for $m(k)$ have been studied by Kurokawa and Ochiai in \cite{KO}, and by Rogers and the author in \cite{LR}. The simplest ones are given as follows:

\begin{thm} We have the following functional equations for $m(k)$:
\begin{itemize}
 \item \cite{KO}: For $h\in\mathbb{R}\backslash \{0\}$:
\begin{equation} \label{eq:ko}
m\left(4h^2\right)+m\left(\frac{4}{h^2}\right) = 2
m\left(2\left(h+\frac{1}{h}\right)\right).
\end{equation}
\item \cite{LR}: If $h \not = 0$, and $|h|<1$:
\begin{equation} \label{eq:lr}
m\left(2\left(h+\frac{1}{h}\right)\right)+m\left(2\left(\ii
h+\frac{1}{\ii h}\right)\right)=  m\left(\frac{4}{h^2}\right).
\end{equation}
\end{itemize}
 \end{thm}
 
If we set $h =\frac{1}{\sqrt{2}}$ in both identities, we obtain
\[m\left(2\right)+m\left(8\right) = 2
m\left(3 \sqrt{2}\right),\]
\[m\left(3 \sqrt{2}\right)+m\left(\ii \sqrt{2}\right)=  m\left(8\right).\]

Similarly, if we set $h=\frac{1}{2}$, we obtain
\[m\left(1\right)+m\left(16\right) = 2
m\left(5\right),\]
\[m\left(5\right)+m\left(-3 \ii \right)=  m\left(16\right).\]
Thus, in order to prove $(\ref{eqA})$ and $(\ref{eqB})$, we need to find one additional equation for each of the above linear systems.

\section{The relationship with the regulator}

In this section, we sometimes write $x_k$ and $y_k$ for $x$ and $y$, so we can distinguish them when we look at different curves. 

After the works of Deninger \cite{D} and Rodriguez-Villegas \cite{RV}, we write
\[m(k) = \frac{1}{2\pi}r_k(\{x_k,y_k\}),\]
were $r_k$ is a period of the regulator in the symbol $\{x_k,y_k\} \in K_2(\mathcal{E}_k)$. For our purposes, we can reduce to $K_2(\C(E_k))$, so that $x_k, y_k$ are elements of $\C(E_k)$. See \cite{D} and \cite{RV} for general details, and \cite{LR} for the specific treatment of this particular example.

In our context, it is enough to take into account that 
\[r_k(\{x_k,y_k\})= \alpha D_k((x_k)\diamond(y_k)),\]
where $\alpha$ is a constant independent of $k$ and $D_k$ is the elliptic dilogarithm in $E_k$ constructed by Bloch (see \cite{Blo}).

We will briefly explain the meaning of $(x) \diamond (y)$. Let $E$ be an elliptic curve with $x,y \in \C(E)$. Consider the divisors
\[(x) = \sum a_S (S), \qquad (y) = \sum b_T (T).\] 
Now define
\[(x)\diamond (y) = \sum a_S b_T (S-T).\]
This is an element in 
\[\Z[E(\C)]^- = \Z[E(\C)]/ \sim ,\]
where the equivalence relation stands for $(-T) \sim -(T)$. 

Thus, the Mahler measure depends just on $D_k$ and $(x_k)\diamond(y_k)$. For example, if the elliptic curves are isomorphic, $D_k$ does not change and the Mahler measure only depends on  $(x_k)\diamond(y_k)$. This idea was discovered by Rodriguez-Villegas \cite{RV2}, and also used by Bertin \cite{Be}. We applied this idea again in \cite{LR}, to isogenous elliptic curves, in order to prove identities like $(\ref{eq:lr})$.

A Weierstrass model for $E_k$ is given by
\[Y^2=X\left(X^2+\left(\frac{k^2}{4}-2\right)X +1\right),\]
where
\[x=\frac{kX-2Y}{2X(X-1)}, \qquad y=\frac{kX+2Y}{2X(X-1)}.\]
It is not hard to see that $E_k(\Q(k))_\mathrm{tor} \cong \Z/4\Z$. To fix notation, we will denote a generator by \[P=\left(1,\frac{k}{2}\right).\] Then we have $2P=(0,0)$. Eventually, we will perform computations in the curve with parameter $k=h + \frac{1}{h}$. In this curve, we will denote  \[Q=\left(-\frac{1}{h^2},0\right),\] which is a point of order 2. Notice that $P+Q=\left(-1, h-\frac{1}{h}\right)$ and $2P+Q=\left(-h^2,0\right)$.

In \cite{LR} we prove
\[(x) \diamond (y) = 8(P).\]

Consider the isomorphism
\[\phi: E_{2\left(h+\frac{1}{h}\right)} \rightarrow E_{2\left(\ii
h+\frac{1}{\ii h}\right)}, \qquad (X,Y) \rightarrow (-X, \ii Y),
 \]
which relates two of the curves in equation (\ref{eq:lr}). We use this isomorphism to pull the rational functions $x, y \in \C\left( E_{2\left(\ii
h+\frac{1}{\ii h}\right)}\right)$ back to $\C\left( E_{2\left(
h+\frac{1}{h}\right)}\right)$:

\[r_{2\left(\ii
h+\frac{1}{\ii h}\right)}(\{x,y\}) = r_{2\left(h+\frac{1}{h}\right)}(\{x \circ \phi ,y \circ \phi\}).\]

On the other hand, it is easy to see that
\[(x \circ \phi) \diamond (y \circ \phi) = 8(P+Q).\]

\section{Relationships between divisors}

From the previous section, the problem reduces to finding relations between $(P)$ and $(P+Q)$ in $\Z\left[E_{2\left(h+\frac{1}{h}\right)}(\C) \right]^-$. In order to do that, we will look for elements that are trivial in $K_2\left(\C\left(E_{2\left(h+\frac{1}{h}\right)}\right) \right)$. In other words, we will find combinations of Steinberg symbols $\{g, 1-g\}$ with $g \in \C\left(E_{2\left(h+\frac{1}{h}\right)}\right)$, such that the corresponding combination $(g)\diamond (1-g)$ yields a linear combination of $(P)$ and $(P+Q)$. Since $\{g, 1-g\}$ is trivial in $K$-theory, we conclude that $(g)\diamond (1-g) \sim 0$, yielding a linear combination involving $(P)$ and $(P+Q)$.

Consider the function
\[f=\frac{Y}{2h}+\left(\frac{1}{2}-\frac{1}{2h^2}\right) X.\]
We have 
\[1-f = 1-\frac{Y}{2h}-\left(\frac{1}{2}-\frac{1}{2h^2}\right) X.\]
Then
\[ (f) = (2P)+2(P+Q) - 3 O, \qquad (1-f)= (P) +(A) +(B) -3 O,\]
where 
\[ A= \left(\frac{-3+\sqrt{9-16h^2}}{2}, \frac{7h}{2}-\frac{3}{2h} - \left(h-\frac{1}{h} \right)\frac{\sqrt{9-16h^2}}{2} \right),\]
\[B = \left(\frac{-3-\sqrt{9-16h^2}}{2}, \frac{7h}{2}-\frac{3}{2h} + \left(h-\frac{1}{h} \right)\frac{\sqrt{9-16h^2}}{2} \right).\]

In particular, for $h=\frac{1}{\sqrt{2}}$, we get
\[ A = 3P+Q, \qquad B = Q,\]
implying
\[(f)\diamond(1-f) = 6(P)-10(P+Q) \sim 0\]
yielding the expected relation.

On the other hand, for $h=\frac{1}{2}$, our function $f$ becomes
\[f=Y-\frac{3}{2} X.\]
%\[1-f=1-Y+\frac{3}{2}X\]

In this case, $A$ and $B$ are given by:
\[A=\left( -\frac{3-\sqrt{5}}{2},-\frac{5-3\sqrt{5}}{4}\right), \qquad B=\left( -\frac{3+\sqrt{5}}{2},-\frac{5+3\sqrt{5}}{4}\right).\] 

In particular, we have the relations
\[2A=2B=P, \quad B-A = 2P, \quad A+B=-P.\]

We obtain
\begin{eqnarray*}
(f) \diamond (1-f) &=& (P) + (2P-A) +(2P-B) -3 (2P)+2(Q) +2(P+Q-A)  \\
&&+2(P+Q-B)-6(P+Q)-3(-P)-3(-A)-3(-B)+9O\\
&=& 2(Q+A) +2(Q+B)-6(P+Q)+4(P)+2(A)+2(B).\\
\end{eqnarray*}

We need further relations among the divisors $(A)$, $(B)$. Thus we consider the following function

\[g=\frac{\sqrt{5}-1}{10}Y+\frac{3+\sqrt{5}}{20}(X+4),\]
 
\[1-g = 1- \frac{\sqrt{5}-1}{10}Y-\frac{3+\sqrt{5}}{20}(X+4).\]

We have 
%\[ (f) = (2P)+2(P+Q) - 3 O\]
%\[(1-f) = (P) + (A) +(B) - 3O\]

\[(g) =(Q) +(A) +(-Q-A) -3O, \qquad (1-g)=(-P) +2(B)-3O.\]

The diamond operation yields a new relation:
\begin{eqnarray*}
(g) \diamond (1-g) &=& (Q+P) +2(Q-B)-3(Q) +(A+P)+ 2(A-B)-3(A)\\
&&+(-Q-A+P)+2(-Q-A-B) -3(-Q-A) -3(P) -6(-B) +9O\\
&=&3(Q+P) -2(Q+B)-3(A)+4(Q+A) -3(P) +5(B).\\
\end{eqnarray*}

In order to get more relations, we apply the Galois conjugate,
\[(g^\sigma) \diamond (1-g^\sigma)= 3(Q+P) -2(Q+A)-3(B)+4(Q+B) -3(P) +5(A).\]

The last two equations yield
\[(g)\diamond (1-g)+  (g^\sigma) \diamond (1-g^\sigma) = 6(Q+P)+2(Q+A) +2(Q+B)+2(A)+2(B) -6(P).\]

Finally, we obtain
\[(f)\diamond (1-f) -(g)\diamond (1-g)-  (g^\sigma) \diamond (1-g^\sigma)= -12(Q+P) +10(P) \sim 0.\]

\section{Conclusion of the proof}

Given a relationship of the form
\[a (P) \sim b (P+Q),\]
we get 
\[ a r_{2\left(h+\frac{1}{h}\right)}\left(\left \{x_{2\left(h+\frac{1}{h}\right)} ,y_{2\left(h+\frac{1}{h}\right)} \right \}\right ) = b r_{2\left(\ii
h+\frac{1}{\ii h}\right)}\left(\left \{x_{2\left(\ii h+\frac{1}{\ii h}\right)} ,y_{2\left(\ii h+\frac{1}{\ii h}\right)} \right \}\right ),\]
and
\[ a m\left(2\left(h+\frac{1}{h}\right)\right) = b m\left( 2\left(\ii
h+\frac{1}{\ii h}\right)\right).\]

Thus, for $h=\frac{1}{\sqrt{2}}$, we recover
\[m(8)= \frac{8}{5} m\left(3\sqrt{2}\right) = \frac{8}{3}  m\left(\ii\sqrt{2}\right) = 4 m(2) .\]
For $h=\frac{1}{2}$, we conclude
\[m(16)= \frac{11}{6} m(5) = \frac{11}{5} m(-3 \ii) = 11 m(1).\]

\[m(5) = 6 m(1).\]

\qed

Questions that remain open are how to predict identities such as (\ref{eqA}) and (\ref{eqB}) and, more precisely, to list all such identities.

\section*{Acknowledgments}
The author would like to thank Herbert Gangl for  his encouragement with this problem. The author is also grateful to David Boyd, Mathew Rogers, and Fernando Rodriguez-Villegas for helpful discussions. Thanks are also due to the referee whose constructive comments have improved the exposition of the paper.

This research was supported by University of Alberta Fac. Sci. Startup Grant N031000610 and NSERC Discovery Grant 355412-2008

%%%%%

%%%%%
\end{document}